\documentclass[a4papper,oneside,onecolumn,12pt]{article}
\usepackage[latin2]{inputenc}
\usepackage{amssymb}
\usepackage{paralist}
\usepackage{amsmath}
\usepackage{amsfonts}
\usepackage{graphicx}
\usepackage{wrapfig}
\usepackage{tikz}
\usepackage{hyperref}

\usetikzlibrary{arrows}
\usetikzlibrary{decorations.pathmorphing}
\usetikzlibrary{snakes}
\usetikzlibrary{positioning}

\newtheorem{prop}{Proposition}[section]
\newtheorem{thm}{Theorem}[section]
\newtheorem{lemma}{Lemma}[section]
\newtheorem{defi}{Definition}[section]

\newtheorem{conj}{Open problem}
\newtheorem{obs}{Observation}[section]

\newcommand{\tab}{\ \ \ \ }

\begin{document}

\vspace{3cm}

\centerline {\Large{\bf Shattering-extremal set systems of $VC$
dimension at most $2$}}

%

\centerline{}

\centerline{\bf {Tam\'as M\'esz\'aros}}

\centerline{}

\centerline{Department of Mathematics, Central European
University}

\centerline{Institute of Mathematics, Budapest University of
Technology and Economics}

\centerline{tmeszaros87@gmail.com}


\centerline{}

\centerline{\bf {Lajos R\'onyai}\footnote{Research supported in
part by OTKA grant NK105645.}}

\centerline{}

\centerline{Computer and Automation Research Institute, Hungarian
Academy of Sciences}

\centerline{Institute of Mathematics, Budapest University of
Technology and Economics}

\centerline{lajos@ilab.sztaki.hu}



\begin{abstract}
We say that a set system $\mathcal{F}\subseteq 2^{[n]}$ shatters a
given set $S\subseteq [n]$ if $2^S=\{F \cap S : F \in
\mathcal{F}\}$. The Sauer inequality states that in general, a set
system $\mathcal{F}$ shatters at least $|\mathcal{F}|$ sets. Here
we concentrate on the case of equality. A set system is called
shattering-extremal if it shatters exactly $|\mathcal{F}|$ sets.
In this paper we characterize shattering-extremal set systems of
Vapnik-Chervonenkis dimension $2$ in terms of their inclusion
graphs, and as a corollary we answer an open question from
\cite{VC1} about leaving out elements from shattering-extremal set
systems in the case of families of Vapnik-Chervonenkis dimension
$2$.
\end{abstract}

\section{Introduction}
\tab Throughout this paper $n$ will be a positive integer, the set
$\{1,2,\dots,n\}$ will be referred to shortly as $[n]$ and the
power set of any set $S\subseteq [n]$ will be denoted by $2^{S}$.
For a set system $\mathcal{F}\subseteq 2^{[n]}$ we will write
$supp(\mathcal{F})$ for its support, i.e.
$supp(\mathcal{F})=\bigcup_{F\in \mathcal{F}} F$.

The central notion of our study is \emph{shattering}.

\begin{defi}
A set system $\mathcal{F}\subseteq 2^{[n]}$ shatters a given set
$S\subseteq [n]$ if
\[2^S=\{F \cap S : F \in \mathcal{F}\}.\]
\end{defi}

The family of subsets of $[n]$ shattered by $\mathcal{F}$ is
denoted by $Sh(\mathcal{F})$. The following inequality states that
in general, a set system $\mathcal{F}$ shatters at least
$|\mathcal{F}|$ sets.

\begin{prop}\label{Sauer inequality}
$|Sh(\mathcal{F})|\geq |\mathcal{F}|$ for every set system
$\mathcal{F}\subseteq 2^{[n]}$.
\end{prop}

The statement was proved by several authors, e.g. Aharoni and
Holzman \cite{Aharoni-Holzman}, Pajor \cite{Pajor}, Sauer
\cite{Sauer}, Shelah \cite{Shelah}. Often it is referred to as the
Sauer inequality. Here we are interested in the case of equality.

\begin{defi}
A set systems $\mathcal{F}\subseteq 2^{[n]}$ is
shattering-extremal, or s-extremal for short, if it shatters
exactly $|\mathcal{F}|$ sets, i.e.
$|\mathcal{F}|=|Sh(\mathcal{F})|$.
\end{defi}

Many interesting results have been obtained in connection with
these combinatorial objects, among others by Bollobás, Leader and
Radcliffe in \cite{Reverse Kleitman}, by Bollobás and Radcliffe in
\cite{Defect sauer}, by Frankl in \cite{Frankl}. Füredi and Quinn
in \cite{Furedi}, and recently Kozma and Moran in \cite{Kozma
Moran} provided interesting examples of s-extremal set systems.
Anstee, Rónyai and Sali in \cite{Shattering news} related
shattering to standard monomials of vanishing ideals, and based on
this relation, the present authors in \cite{diploma} and in
\cite{Springer} developed algebraic methods for the investigation
of s-extremal families.

\begin{defi}
The inclusion graph of a set system $\mathcal{F}\subseteq
2^{[n]}$, denoted by $\mathbb{G}_{\mathcal{F}}$, is the simple
directed edge labelled graph whose vertices are the elements of
$\mathcal{F}$, and there is a directed edge with label $j\in [n]$
going from $G$ to $F$ exactly when $F=G\cup \{j\}$.
\end{defi}

$\mathbb{G}_{\mathcal{F}}$ is actually the Hasse diagram of the
poset $\mathcal{F}$ with edges directed and labelled in a natural
way. The inclusion graph of the complete set system $2^{[n]}$ will
be denoted by $\mathbb{H}_n$. The undirected version of
$\mathbb{H}_n$ is often referred to as the Hamming graph $H(n,2)$,
or as the hypercube of dimension $n$, whose vertices are all $0-1$
vectors of length $n$, and two vertices are adjacent iff they
differ in exactly one coordinate. When computing distances between
vertices in the inclusion graph $\mathbb{G}_{\mathcal{F}}$ we
forget about the direction of edges, and define the distance
between vertices $F,G\in \mathcal{F}$, denoted by
$d_{\mathbb{G}_{\mathcal{F}}}(F,G)$, as their graph distance in
the undirected version of $\mathbb{G}_{\mathcal{F}}$, i.e. the
length of the shortest path between them in the undirected version
of $\mathbb{G}_{\mathcal{F}}$. Similarly, some edges in
$\mathbb{G}_\mathcal{F}$ form a path between two vertices if they
do so in the undirected version of $\mathbb{G}_\mathcal{F}$. For
example, the distance between two vertices $F,G\subseteq [n]$ in
$\mathbb{H}_n$ is just the size of the symmetric difference
$F\bigtriangleup G$, i.e. $d_{\mathbb{H}_n}(F,G)=|F\bigtriangleup
G|$. As a consequence, when only distances of vertices will be
considered, and the context will allow, we omit the directions of
edges to avoid unnecessary case analysis, and will specify edges
by merely listing their endpoints.

\begin{defi}
The Vapnik-Chervonenkis dimension of a set system
$\mathcal{F}\subseteq 2^{[n]}$, denoted by
$dim_{VC}(\mathcal{F})$, is the maximum cardinality of a set
shattered by $\mathcal{F}$.
\end{defi}

The general task of giving a good description of s-extremal
systems seems to be too complex at this point, therefore we
restrict our attention to the simplest cases, where the
$VC$-dimension of $\mathcal{F}$ is small. S-extremal systems,
where the $VC$-dimension is at most $1$ were fully described in
\cite{VC1}.

\begin{prop}\label{extr char VC-dim=1} (See \cite{VC1}.)
A set system $\mathcal{F}\subseteq 2^{[n]}$ is s-extremal and of
$VC$ dimension at most $1$ iff $G_{\mathcal{F}}$ is a tree and all
labels on the edges are different.
\end{prop}

Proposition \ref{extr char VC-dim=1} can also be interpreted as
follows:

\begin{prop}\label{VC1 correspondance} (See \cite{VC1})
There is a one-to-one correspondence between s-extremal families
$\mathcal{F}\subseteq 2^{[n]}$ of $VC$-dimension $1$ with
$supp(\mathcal{F})=[n]$, $\cap_{F\in \mathcal{F}} F=\emptyset$ and
directed edge-labelled trees on $n+1$ vertices, all edges having a
different label from $[n]$.
\end{prop}

Note that the assumptions $supp(\mathcal{F})=[n]$ and $\cap_{F\in
\mathcal{F}} F=\emptyset$ are not restrictive. Both of them can be
assumed to hold without loss of generality, otherwise one could
omit common elements and then restrict the ground set to
$supp(\mathcal{F})$.

In this paper we continue the work initiated in \cite{VC1}, and
characterize s-extremal set systems of $VC$-dimension at most $2$.
We do this by providing an algorithmic procedure for constructing
the inclusion graphs of all such set systems. This
characterization then allows us to answer an open question, posed
in \cite{VC1}, about leaving out elements from such set systems.

The paper is organized as follows. After the introduction in
Section \ref{preliminaries} we investigate the properties of
shattering and its connection to inclusion graphs. Next, in
Section \ref{section constr} we propose a building process for
extremal families and investigate its properties. Based on this
building process in Section \ref{section main results} we present
and prove our main results. Finally in Section \ref{section
remarks} we make some concluding remarks concerning future work.

\section{Preliminaries}\label{preliminaries}

To start with, we first introduce a useful subdivision of set
systems.

\begin{defi}\label{stand subdiv def}
The standard subdivision of a set system $\mathcal{F}\subseteq
2^{[n]}$ with respect to an element $i\in [n]$ consists of the
following two set systems:
\begin{center}
$\mathcal{F}_0=\{F : F \in \mathcal{F} ; i \notin F\}\subseteq 2^{[n]\backslash\{i\}}$ and\\
$\mathcal{F}_1=\{F\backslash \{i\} : F \in \mathcal{F} ; i \in
F\}\subseteq 2^{[n]\backslash\{i\}}$.
\end{center}
\end{defi}

\noindent For the sake of completeness we provide a possible proof
of Proposition \ref{Sauer inequality}, whose main idea will be
useful later on.\\

\textbf{Proof:}(of Proposition \ref{Sauer inequality}) We will
prove this statement by induction on $n$. For $n=1$ the statement
is trivial. Now suppose that $n>1$, and consider the standard
subdivision of $\mathcal{F}$ with respect to the element $n$. Note
that $\mathcal{F}_0,\mathcal{F}_1\subseteq 2^{[n-1]}$ and hence by
the induction hypothesis we have $|Sh(\mathcal{F}_0)|\geq
|\mathcal{F}_0|$ and $|Sh(\mathcal{F}_1)|\geq |\mathcal{F}_1|$.
Moreover $|\mathcal{F}|=|\mathcal{F}_0|+|\mathcal{F}_1|$,
$Sh(\mathcal{F}_0)\cup Sh(\mathcal{F}_1) \subseteq
Sh(\mathcal{F})$ and if $S \in Sh(\mathcal{F}_0)\cap
Sh(\mathcal{F}_1)$, then according to the definition of
$\mathcal{F}_0$ and $\mathcal{F}_1$ we have $S\cup \{n\} \in
Sh(\mathcal{F})$. Summarizing
\[|Sh(\mathcal{F})|\geq
|Sh(\mathcal{F}_0)|+|Sh(\mathcal{F}_1)|\geq
|\mathcal{F}_0|+|\mathcal{F}_1|=|\mathcal{F}|. \blacksquare\]

From the proof of Proposition \ref{Sauer inequality} it is easy to
see, that if $\mathcal{F}$ is s-extremal, then so are the systems
$\mathcal{F}_0$ and $\mathcal{F}_1$ in the standard subdivision
with respect to any element $i\in [n]$. Iterating this for an
s-extremal system $\mathcal{F}\subseteq 2^{[n]}$ we get that for
all pairs of sets $A\subseteq B\subseteq [n]$, the system

\[\{F\backslash A\ |\ F\in \mathcal{F},\ A\subseteq F\subseteq B\}\]

\noindent is s-extremal. Moreover if in the above system we add
$A$ to every set, then the family of shattered sets remains
unchanged, hence we get that the subsystem

\[\mathcal{F}_{A,B}=\{F\ |\ F\in \mathcal{F},\ A\subseteq F\subseteq B\}\subseteq \mathcal{F}\]

\noindent is also s-extremal.

In \cite{Reverse Kleitman} and \cite{Defect sauer} a different
version of shattering, \emph{strong shattering} is introduced .

\begin{defi}
A set system $\mathcal{F}\subseteq2^{[n]}$ strongly shatters the
set $F\subseteq [n]$, if there exists $I\subseteq [n]\backslash F$
such that

\[2^F+I=\{H\cup I\ |\ H\subseteq F\}\subseteq \mathcal{F}.\]
\end{defi}

The family of all sets strongly shattered by some set system
$\mathcal{F}$ is denoted by $st(\mathcal{F})$. Clearly
$st(\mathcal{F})\subseteq Sh(\mathcal{F})$, both $Sh(\mathcal{F})$
and $st(\mathcal{F})$ are down sets and both families are
monotone, meaning that if $\mathcal{F}\subseteq \mathcal{F}'$ are
set systems then $Sh(\mathcal{F})\subseteq Sh(\mathcal{F'})$ and
$st(\mathcal{F})\subseteq st(\mathcal{F'})$. For the size of
$st(\mathcal{F})$ one can prove the so called reverse Sauer
inequality:

\begin{prop}
(see \cite{Reverse Kleitman}) $|st(\mathcal{F})|\leq
|\mathcal{F}|$ for every set system $\mathcal{F}\subseteq
2^{[n]}$.
\end{prop}

Bollobás and Radcliffe in \cite{Defect sauer} obtained several
important results concerning shattering and strong shattering,
including:

\begin{prop}
(see \cite{Defect sauer}, Theorem $2$) $\mathcal{F}\subseteq
2^{[n]}$ is extremal with respect to the Sauer inequality (i.e. is
shattering-extremal) iff it is extremal with respect to the
reverse Sauer inequality i.e. $|st(\mathcal{F})|=|\mathcal{F}|\
\Longleftrightarrow \ |Sh(\mathcal{F})|=|\mathcal{F}|$.
\end{prop}

Since the two extremal cases coincide, we will call such set
systems shortly just extremal. As a consequence of the above
facts, we get, that for extremal systems we have
$st(\mathcal{F})=Sh(\mathcal{F})$.

For $i\in [n]$ let $\varphi_i$ be the $i$th \emph{bit flip
operation}, i.e. for $F\in 2^{[n]}$ we have

\[\varphi_i(F)=F\triangle\{i\}=\left\{
\begin{array}{ll}
  \mbox{$F\backslash\{i\}$} & \mbox{ if  $i\in F$} \\
  \mbox{$F\cup{\{i\}}$} & \mbox{ if $i\notin F$ }\\
\end{array}
\right.\]

\noindent and $\varphi_i(\mathcal{F})=\{\varphi_i(F)\ |\ F\in
\mathcal{F}\}$. The family of shattered sets is trivially
invariant under the bit flip operation, i.e.
$Sh(\mathcal{F})=Sh(\varphi_i(\mathcal{F}))$ for all $i\in [n]$,
and hence so is extremality. This means that when dealing with a
nonempty set system $\mathcal{F}$, and examining its extremality,
we can assume that $\emptyset\in \mathcal{F}$, otherwise we could
apply bit flips to it, to bring $\emptyset$ inside.

In terms of the inclusion graph, $\varphi_i$ flips the directions
of edges with label $i$, i.e. there is a bijection between the
vertices of $\mathbb{G}_{\mathcal{F}}$ and
$\mathbb{G}_{\varphi_i(\mathcal{F})}$ that preserves all edges
with label different from $i$, and reverses edges with label $i$.
This bijection is simply given by the reflection with respect to
the hyperplane $x_i=\frac{1}{2}$ in the Hamming graph, when viewed
as a subset of $\mathbb{R}^n$.

Note that for any set system $\mathcal{F}\subseteq 2^{[n]}$, the
identity map naturally embeds the inclusion graph
$\mathbb{G}_{\mathcal{F}}$ into $\mathbb{H}_n$. We say that the
inclusion graph $\mathbb{G}_{\mathcal{F}}$ is isometrically
embedded (into $\mathbb{H}_n$), if this embedding is an isometry,
meaning that for arbitrary $F,G\in \mathcal{F}$ we have
$d_{\mathbb{G}_\mathcal{F}}(F,G)=d_{\mathbb{H}_n}(F,G)$, i.e.
there is a path of length $d_{\mathbb{H}_n}(F,G)=|F\bigtriangleup
G|$ between $F$ and $G$ inside the undirected version of
$\mathbb{G}_{\mathcal{F}}$. Greco in \cite{greco} proved the
following:

\begin{prop}\label{extr->isom}
If $\mathcal{F} \subseteq 2^{[n]}$ is extremal, then
$\mathbb{G}_{\mathcal{F}}$ is isometrically embedded.
\end{prop}

As this fact will be used several times, we provide the reader
with a simple
proof from \cite{diploma}:\\

\textbf{Proof:} Suppose the contrary, namely that
$\mathbb{G}_{\mathcal{F}}$ is not isometrically embedded. Then
there exist sets $A,B\in \mathcal{F}$ such that
$d_{\mathbb{H}_n}(A,B)=k < d_{\mathbb{G}_{\mathcal{F}}}(A,B)$.
Suppose that $A$ and $B$ are such that $k$ is minimal. Clearly
$k\geq 2$. W.l.o.g we may suppose that $A=\emptyset$ and $|B|=k$,
otherwise one could apply bit flips to the set system to achieve
this. Note that distances both in $\mathbb{G}_{\mathcal{F}}$ and
in $\mathbb{H}_n$ are invariant under bit flips.

We claim that there is no set $C \in \mathcal{F}$ different from
$A$ with $C \varsubsetneq B$. Indeed suppose such $C$ exists, then
\[d_{\mathbb{H}_n}(A,C)+d_{\mathbb{H}_n}(C,B)=d_{\mathbb{H}_n}(A,B)=k
< d_{\mathbb{G}_{\mathcal{F}}}(A,B)\leq
d_{\mathbb{G}_{\mathcal{F}}}(A,C)+d_{\mathbb{G}_{\mathcal{F}}}(C,B).\]
From this we have either $d_{\mathbb{H}_n}(A,C)<
d_{\mathbb{G}_{\mathcal{F}}}(A,C)$ or $d_{\mathbb{H}_n}(C,B)<
d_{\mathbb{G}_{\mathcal{F}}}(C,B)$. Since
$d_{\mathbb{H}_n}(A,C),d_{\mathbb{H}_n}(C,B)<k$ we get a
contradiction in both cases with the minimality of $k$.

Now since $\mathcal{F}$ is extremal, so must be
$\mathcal{F}_{\emptyset,B}$. However in our case
$\mathcal{F}_{\emptyset,B}=\{\emptyset, B\}$, and so if
$B=\{b_1,\dots,b_k\}$, then
$Sh(\mathcal{F}_{\emptyset,B})=\{\emptyset,\{b_1\},\dots,\{b_k\}\}$.
Counting cardinalities we get that
$|Sh(\mathcal{F}_{\emptyset,B})|=|B|+1=k+1\geq 3 >
2=|\mathcal{F}_{\emptyset,B}|$, implying that
$\mathcal{F}_{\emptyset,B}$ cannot be extremal. This contradiction
finishes the proof.
$\blacksquare$\\

It is easy to see that $S\in st(\mathcal{F})$ (and so in the
extremal case $S\in Sh(\mathcal{F})$) is just equivalent to the
fact that $\mathbb{G}_{2^S}$ is isomorphic to a subgraph of
$\mathbb{G}_{\mathcal{F}}$ as a directed edge labelled graph, i.e.
there exists a bijection between the vertices of
$\mathbb{G}_{2^S}$ and $2^{|S|}$ vertices of
$\mathbb{G}_{\mathcal{F}}$ preserving edges, edge labels and edge
directions. If this happens, we will say, that there is a copy of
$\mathbb{G}_{2^S}$ in $\mathbb{G}_{\mathcal{F}}$.

Suppose that for a set $S\subseteq [n]$ there are $2$ different
copies of $\mathbb{G}_{2^S}$ in $\mathbb{G}_{\mathcal{F}}$, i.e.
there are two different sets $I_1,I_2\subseteq [n]\backslash S$
such that $2^{S}+I_1,2^{S}+I_2\subseteq \mathcal{F}$. Since
$I_1\neq I_2$, there must be an element $\alpha\notin S$ such that
$\alpha\in I_1\triangle I_2$. For this element $\alpha$ we clearly
have that $\mathcal{F}$ shatters $S\cup\{\alpha\}$.

\begin{obs}\label{max el uniq st}
If $\mathcal{F}\subseteq 2^{[n]}$ is extremal and $S\subseteq [n]$
is a maximal element in $st(\mathcal{F})=Sh(\mathcal{F})$, in the
sense that $S\in st(\mathcal{F})=Sh(\mathcal{F})$ and for all
$S'\varsupsetneq S$ we have $S'\notin
st(\mathcal{F})=Sh(\mathcal{F})$, then $S$ is uniquely strongly
shattered, i.e. there is one unique copy of $\mathbb{G}_{2^S}$ in
$\mathbb{G}_{\mathcal{F}}$
\end{obs}

Indeed, by the earlier reasoning, multiple copies would result a
contradiction with the maximality of $S$.

\section{Construction of extremal families}\label{section constr}

In this section we will describe and study a process for building
up an extremal set system on the ground set $[n]$ together with
its inclusion graph. First we describe the building process for
the set system and then study how the inclusion graph evolves in
the meantime. Let Step $0$ be the initialization, after which we
are given the set system $\{\emptyset\}$. Now suppose we are given
a set system $\mathcal{F}$ and consider the following two types of
operations to enlarge $\mathcal{F}$:

\begin{itemize}
\item \textbf{Step A} - If such exists, take an element $\alpha\in
[n]\backslash supp(\mathcal{F})$ together with a set $W\in
\mathcal{F}$ and add the set $V=\{W,\alpha\}$ to $\mathcal{F}$.

Note that the singleton $\{\alpha\}$ is strongly shattered by
$\mathcal{F}\cup \{V\}$, as shown by the sets $W$ and $V$, but is
not by $\mathcal{F}$, as by assumption $\alpha\notin
supp(\mathcal{F})$.

\item \textbf{Step B} - If there exist, take two elements
$\alpha,\beta\in supp(\mathcal{F})$ such that
$\{\alpha,\beta\}\notin st(\mathcal{F})$, together with sets
$P,W,Q\in \mathcal{F}$ such that $Q\bigtriangleup W=\{\alpha\}$
and $P\bigtriangleup W=\{\beta\}$. Let $V=W\bigtriangleup
\{\alpha,\beta\}$. $V$ is also the unique set satisfying
$P\bigtriangleup V=\{\alpha\}$ and $Q \bigtriangleup V=\{\beta\}$.
For these sets we have that $\{P,W,Q,V\}=W\cap
V+2^{\{\alpha,\beta\}}=P\cap Q+2^{\{\alpha,\beta\}}$, and hence
$V$ cannot belong to $\mathcal{F}$, otherwise the sets $P,W,Q,V$
would strongly shatter $\{\alpha,\beta\}$, contradicting our
assumption. Therefore, it is reasonable to add $V$ to
$\mathcal{F}$.

Note that the set $\{\alpha,\beta\}$ is strongly shattered by
$\mathcal{F}\cup \{V\}$, as shown by the sets $P,W,Q$ and $V$, but
is not by $\mathcal{F}$ by assumption.

\end{itemize}

Let $\mathcal{E}$ be the collection of all set systems
$\mathcal{F}$ that can be built up starting with Step $0$ and then
using steps of type A and B in an arbitrary but valid order.

\begin{lemma}\label{steps keep extremality}
Any set system $\mathcal{F}\in\mathcal{E}$ is extremal and
$dim_{VC}(\mathcal{F})\leq 2$.
\end{lemma}
\textbf{Proof:} We will use induction on the size of
$\mathcal{F}$. If $|\mathcal{F}|=1$ then necessarily
$\mathcal{F}=\{\emptyset\}$, which is clearly extremal and
$dim_{VC}(\mathcal{F})=0$. Now suppose we know the result for all
members of $\mathcal{E}$ of size at most $m\geq 1$, and consider a
system $\mathcal{F}\in \mathcal{E}$ of size $m+1$. As
$\mathcal{F}\in \mathcal{E}$ it can be built up starting from
$\{\emptyset\}$ using Steps A and B. Fix one such building
process, and let $\mathcal{F'}$ be the set system before the last
building step. As noted previously, independently of the type of
the last step there is a set $S$ that is strongly shattered by
$\mathcal{F}$ but is not strongly shattered by $\mathcal{F}'$.
$S$ is either a singleton or a set of size $2$,
depending on the type of the last step. By the induction
hypothesis $\mathcal{F}'$ is extremal and
$dim_{VC}(\mathcal{F}')\leq 2$. Using the reverse Sauer
inequality we get that
\[|\mathcal{F}'|=|st(\mathcal{F}')|<|st(\mathcal{F})|\leq |\mathcal{F}|=|\mathcal{F}'|+1,\]
what is possible only if $|st(\mathcal{F})|=|st(\mathcal{F}')|+1$
and $|st(\mathcal{F})|=|\mathcal{F}|$, in particular $\mathcal{F}$
is extremal.

However in the extremal case the family of shattered sets is the
same as the family of strongly shattered sets, and so the above
reasoning also gives that there is exactly one set that is
shattered by $\mathcal{F}$ and is not shattered by $\mathcal{F}'$,
namely $S$, and so $dim_{VC}(\mathcal{F})\leq
max(dim_{VC}(\mathcal{F}'),|S|)\leq 2$. $\blacksquare$\\

The proof of Lemma \ref{steps keep extremality} also describes how
the family of shattered/strongly shattered sets grows during a
building process. After each step it grows by exactly one new set,
namely by $\{\alpha\}$, if the step considered was Step A with the
label $\alpha$, and by $\{\alpha,\beta\}$, if the step considered
was Step B with labels $\alpha,\beta$. By our assumptions on the
steps it also follows that a valid building process for a set
system $\mathcal{F}\in \mathcal{E}$ cannot involve twice Step A
with the same label $\alpha$, neither twice Step B with the same
pair of labels $\alpha,\beta$, and we also have that
\[Sh(\mathcal{F})=st(\mathcal{F})=\begin{array}{c}
  \Big\{\emptyset\Big\}\bigcup \Big\{\{\alpha\}\ |\ \mbox{Step A is used with label }\alpha\Big\}\bigcup \\
  \Big\{\{\alpha,\beta\}\ |\ \mbox{Step B is used with labels }\alpha\mbox{ and }\beta\Big\} \\
\end{array}.\]

Now consider a valid building process from $\mathcal{E}$, and let us examine,
how the inclusion graph evolves. We use the notation
from the definitions of Steps A and B. Suppose we have already
built up a set system $\mathcal{F}$, and we are given its
inclusion graph $G_{\mathcal{F}}$.

In Step A we add a new vertex, namely $V$ to $G_{\mathcal{F}}$,
together with one new directed edge with label $\alpha$ going from
$W$ to $V$. As $\alpha\notin supp(\mathcal{F})$, $V$ has no other
neighbors in $G_{\mathcal{F}}$. Figure \ref{fig step A} shows Step
A in terms of the inclusion graph.

\begin{figure}
\centering
\begin{tikzpicture}[-,>=stealth',shorten >=1pt,auto,node distance=3cm,
      thick,main node/.style={circle,fill=none,draw,font=\sffamily\tiny\bfseries}]

    \def \radius {2cm}
    \def \margin {8}

    \node[main node, draw, circle] (1) at ({0}:\radius) {$W$};
    \node[main node, draw, circle,red] (2) at ({0}:5) {$V$};

    \draw[loosely dashed] ({5+\margin}:\radius)
         arc ({6+\margin}:{360-\margin}:\radius);

    \draw[->,red] (1) -- (2) node[pos=.5,above] {$\alpha$};

    \end{tikzpicture}
    \caption{Step A}
    \label{fig step A}
\end{figure}

\begin{figure}
\centering
\begin{tikzpicture}[-,>=stealth',shorten >=1pt,auto,node distance=3cm,
      thick,main node/.style={circle,fill=none,draw,font=\sffamily\tiny\bfseries}]

    \def \radius {3cm}
    \def \margin {8}

    \node[main node, draw, circle] (1) at ({45}:\radius) {$P$};
    \node[main node, draw, circle] (2) at ({315}:\radius) {$Q$};
    \node[main node, draw, circle] (3) at ({0}:0) {$W$};
    \node[main node, draw, circle,red] (4) at ({0}:4.2) {$V$};
    \node[main node, draw, circle] (5) at ({135}:1.8) {$Y$};
    \node[main node, draw, circle] (6) at ({225}:1.8) {$X$};
    \node[main node, draw, circle] (7) at ({280}:2.2) {$Z$};

    \draw[loosely dashed] ({45+\margin}:\radius)
         arc ({45+\margin}:{315-\margin}:\radius);

    \draw[-] (1) -- (3) node[pos=.3,left] {$\beta$};
    \draw[-,red] (4) -- (2) node[pos=.5,right] {$\beta$};
    \draw[-,red] (1) -- (4) node[pos=.5,right] {$\alpha$};
    \draw[-] (3) -- (2) node[pos=.8,left] {$\alpha$};
    \draw[-,red] (6) -- (4) node[pos=.7,above] {$\gamma$};
    \draw[-] (6) -- (5) node[pos=.5,left] {$\alpha$};
    \draw[-] (5) -- (1) node[pos=.5,above] {$\gamma$};
    \draw[-] (7) -- (6) node[pos=.5,below] {$\beta$};
    \draw[-] (2) -- (7) node[pos=.6,below] {$\gamma$};

    \end{tikzpicture}
    \caption{Case of Step B}
    \label{fig case step B}
\end{figure}
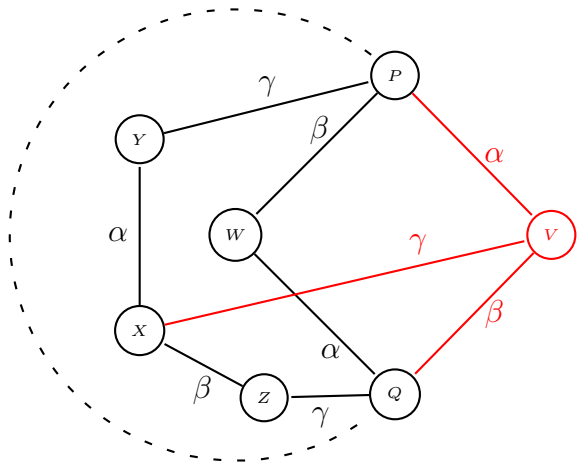

In Step B we also add one new vertex to $G_{\mathcal{F}}$, namely
$V$. As the distance of $V$ from both $P$ and $Q$ is $1$, and
$P\bigtriangleup V=\{\alpha\}$ and $Q \bigtriangleup V=\{\beta\}$,
we have to add at least $2$ new edges, one between $P$ and $V$
with label $\alpha$ and one between $Q$ and $V$ with label
$\beta$. The direction of these edges is predetermined by the
vertices $P,W$ and $Q$. Figure \ref{fig step B} shows all possible
cases for the directions of these edges. We claim that no other edges need to be added, i.e. $V$ has no
other neighbors in $G_{\mathcal{F}}$. Indeed suppose that the
new vertex $V$ has another neighbor $X$ in $G_{\mathcal{F}}$,
different from $P$ and $Q$, that should be connected to it with
some label $\gamma$ different from $\alpha$ and $\beta$. See
Figure \ref{fig case step B}, where edge directions are ignored,
only edge labels are shown.

\begin{figure}
\centering
\begin{tikzpicture}[-,>=stealth',shorten >=1pt,auto,node distance=2cm,
      thick,main node/.style={circle,fill=none,draw,font=\sffamily\tiny\bfseries}]

    \def \radius {2cm}
    \def \margin {8}

    \node[main node] (1) {$W$};
    \node[main node] (2) [above right of=1] {$P$};
    \node[main node] (3) [below right of=1] {$Q$};
    \node[main node] (4) [below right of=2] {$V$};

    \draw[loosely dashed] (2)+(-0.2,0.2) arc (50:310:\radius);

    \draw[->] (1) -- (2) node[pos=.7,left] {$\beta$};
    \draw[<-,red] (4) -- (3) node[pos=.7,right] {$\beta$};
    \draw[<-,red] (2) -- (4) node[pos=.5,right] {$\alpha$};
    \draw[->] (3) -- (1) node[pos=.3,left] {$\alpha$};

    \node[main node] (5) [right=6 cm of 1] {$W$};
    \node[main node] (6) [above right of=5] {$P$};
    \node[main node] (7) [below right of=5] {$Q$};
    \node[main node] (8) [below right of=6] {$V$};

    \draw[loosely dashed] (6)+(-0.2,0.2) arc (50:310:\radius);

    \draw[<-] (5) -- (6) node[pos=.7,left] {$\beta$};
    \draw[->,red] (8) -- (7) node[pos=.7,right] {$\beta$};
    \draw[<-,red] (6) -- (8) node[pos=.3,right] {$\alpha$};
    \draw[->] (7) -- (5) node[pos=.5,left] {$\alpha$};

    \node[main node] (9) [below=5 cm of 1] {$W$};
    \node[main node] (10) [above right of=9] {$P$};
    \node[main node] (11) [below right of=9] {$Q$};
    \node[main node] (12) [below right of=10] {$V$};

    \draw[loosely dashed] (10)+(-0.2,0.2) arc (50:310:\radius);

    \draw[<-] (9) -- (10) node[pos=.7,left] {$\beta$};
    \draw[->,red] (12) -- (11) node[pos=.7,right] {$\beta$};
    \draw[->,red] (10) -- (12) node[pos=.3,right] {$\alpha$};
    \draw[<-] (11) -- (9) node[pos=.5,left] {$\alpha$};

    \node[main node] (13) [right=6 cm of 9] {$W$};
    \node[main node] (14) [above right of=13] {$P$};
    \node[main node] (15) [below right of=13] {$Q$};
    \node[main node] (16) [below right of=14] {$V$};

    \draw[loosely dashed] (14)+(-0.2,0.2) arc (50:310:\radius);

    \draw[->] (13) -- (14) node[pos=.7,left] {$\beta$};
    \draw[<-,red] (16) -- (15) node[pos=.7,right] {$\beta$};
    \draw[->,red] (14) -- (16) node[pos=.3,right] {$\alpha$};
    \draw[<-] (15) -- (13) node[pos=.5,left] {$\alpha$};

    \end{tikzpicture}
    \caption{Step B}
    \label{fig step B}
\end{figure}

Here $d_{\mathbb{H}_n}(P,X)=|P\bigtriangleup
X|=|\{\alpha,\gamma\}|=2$. On the other hand as $\mathcal{F}$ was
built using Steps A and B starting from $\{\emptyset\}$, it is a
member of $\mathcal{E}$, and so by Lemma \ref{build extr syst} it
is extremal. According to Proposition \ref{extr->isom} this
implies that $\mathcal{G}_{\mathcal{F}}$ is isometrically
embedded. This means that there should be a vertex $Y$ in
$\mathbb{G}_{\mathcal{F}}$ connected to both $P$ and $X$ with
edges with labels $\gamma$ and $\alpha$ respectively. The same
reasoning applies for $Q$ and $V$ with some intermediate vertex
$Z$ and edge labels $\beta$, $\gamma$. However in this case,
independently of the directions of the edges, we have
$\{X\cap\{\alpha,\beta\},Y\cap\{\alpha,\beta\},Z\cap\{\alpha,\beta\},W\cap\{\alpha,\beta\}\}=2^{\{\alpha,\beta\}}$,
i.e. the sets $X,Y,Z,W$ shatter the set $\{\alpha,\beta\}$, and so
by the extremality of $\mathcal{F}$ we have that
$\{\alpha,\beta\}$ is also strongly shattered, what contradicts
the assumptions of Step B.

From now on it will depend on the context whether we regard Steps
A and B as building steps for extremal set systems of $VC$
dimension at most $2$ or as building steps for their inclusion
graphs.

\begin{figure}
\centering
\begin{tikzpicture}[-,>=stealth',shorten >=1pt,auto,node distance=2cm,
      thick,main node/.style={rectangle,fill=none,draw,font=\sffamily\tiny\bfseries}]

    \node[main node,red] (1) {$\emptyset$};

    \node[main node] (2) [right=3cm of 1] {$\emptyset$};
    \node[main node, red] (3) [right of=2] {$\{1\}$};

    \node[main node] (4) [right=3 cm of 3] {$\emptyset$};
    \node[main node] (5) [above right of=4] {$\{1\}$};
    \node[main node, red] (6) [below right of=4] {$\{2\}$};

    \node[main node] (7) at (1,-6) {$\emptyset$};
    \node[main node] (8) [above right of=7] {$\{1\}$};
    \node[main node] (9) [below right of=7] {$\{2\}$};
    \node[main node, red] (10) [above right of=9] {$\{2,3\}$};

    \node[main node] (11) [right=6.5 of 7] {$\emptyset$};
    \node[main node, red] (15) [above right of=11] {$\{3\}$};
    \node[main node] (12) [above of=15] {$\{1\}$};
    \node[main node] (13) [below right of=11] {$\{2\}$};
    \node[main node] (14) [above right of=13] {$\{2,3\}$};

    \draw[->,snake] (-2,0) -- (-0.5,0) node[pos=.5,above] {Step $0$};

    \draw[->,snake] (1,0) -- (2.5,0) node[pos=.5,above] {$\begin{array}{c}
      \mbox{Step }A \\
      \mbox{with label }1 \\
    \end{array}$};

    \draw[->,red] (2) -- (3) node[pos=.5,above] {$1$};

    \draw[->,snake] (6.5,0) -- (8,0) node[pos=.5,above] {$\begin{array}{c}
      \mbox{Step }A \\
      \mbox{with label }2 \\
    \end{array}$};

    \draw[->] (4) -- (5) node[pos=.5,above] {$1$};
    \draw[<-,red] (6) -- (4) node[pos=.5,above] {$2$};

    \draw[->,snake] (-2,-6) -- (0,-6) node[pos=.5,above] {$\begin{array}{c}
      \mbox{Step }A \\
      \mbox{with label }3 \\
    \end{array}$};

    \draw[->] (7) -- (8) node[pos=.5,above] {$1$};
    \draw[<-] (9) -- (7) node[pos=.5,above] {$2$};
    \draw[<-,red] (10) -- (9) node[pos=.5,above] {$3$};

    \draw[->,snake] (5,-6) -- (7,-6) node[pos=.5,above] {$\begin{array}{c}
      \mbox{Step }B \\
      \mbox{with labels }2,3 \\
    \end{array}$};

    \draw[->] (11) -- (12) node[pos=.5,above] {$1$};
    \draw[<-] (13) -- (11) node[pos=.5,above] {$2$};
    \draw[<-] (14) -- (13) node[pos=.5,above] {$3$};
    \draw[<-,red] (15) -- (11) node[pos=.5,above] {$3$};
    \draw[->,red] (15) -- (14) node[pos=.5,above] {$2$};
\end{tikzpicture}

    \caption{Example of a building process in $\mathcal{E}$}
    \label{fig example E}
\end{figure}
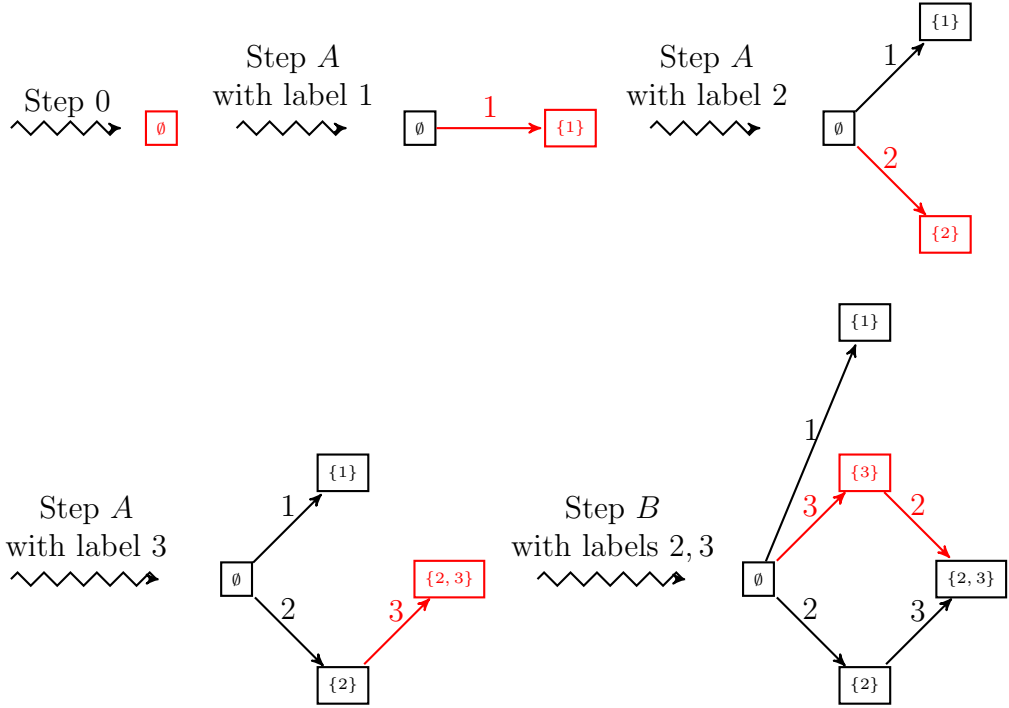

Figure \ref{fig example E} shows a possible building process in
$\mathcal{E}$ for the set system
\[\mathcal{F}=\{\emptyset,\{1\},\{2\},\{3\},\{2,3\}\}\]
in terms of the inclusion graph.

\begin{figure}
\centering
\begin{tikzpicture}[-,>=stealth',shorten >=1pt,auto,node distance=3cm,
      thick,main node/.style={circle,fill=none,draw,font=\sffamily\tiny\bfseries}]

      \node[main node] (1) {$X_1$};
      \node[main node] (2) [below of=1] {$Y_1$};
      \node[main node] (3) [below right of=1] {};
      \node[main node] (4) [below of=3] {};
      \node[main node] (5) [above right of=3] {};
      \node[main node] (6) [below of=5] {};
      \node[main node] (7) [below right of=5] {};
      \node[main node] (8) [below of=7] {};
      \node[main node] (9) [above right of=7] {};
      \node[main node] (10) [below of=9] {};
      \node[main node] (11) [below right of=9] {$X_2$};
      \node[main node] (12) [below of=11] {$Y_2$};

      \draw[color=red] (1) -- (2) node[pos=.5,left] {$\alpha$};
      \draw[] (3) -- (4) node[pos=.5,right] {$\alpha$};
      \draw[] (5) -- (6) node[pos=.5,right] {$\alpha$};
      \draw[] (7) -- (8) node[pos=.5,right] {$\alpha$};
      \draw[] (9) -- (10) node[pos=.5,right] {$\alpha$};
      \draw[color=red] (11) -- (12) node[pos=.5,right] {$\alpha$};

      \draw[] (1) -- (3) node[pos=.5,above] {$\beta_1$};
      \draw[] (3) -- (5) node[pos=.5,above] {$\beta_2$};
      \draw[] (7) -- (9) node[pos=.5,left] {$\beta_{\ell-1}$};
      \draw[] (9) -- (11) node[pos=.5,above] {$\beta_{\ell}$};

      \draw[] (4) -- (2) node[pos=.5,below] {$\beta_1$};
      \draw[] (6) -- (4) node[pos=.5,below] {$\beta_2$};
      \draw[] (10) -- (8) node[pos=.5,right] {$\beta_{\ell-1}$};
      \draw[] (12) -- (10) node[pos=.5,below] {$\beta_{\ell}$};

      \draw[loosely dashed] (5) -- (7);
      \draw[loosely dashed] (6) -- (8);
    \end{tikzpicture}
    \caption{Path of $4$ cycles}
    \label{fig path 4-cycles}
\end{figure}
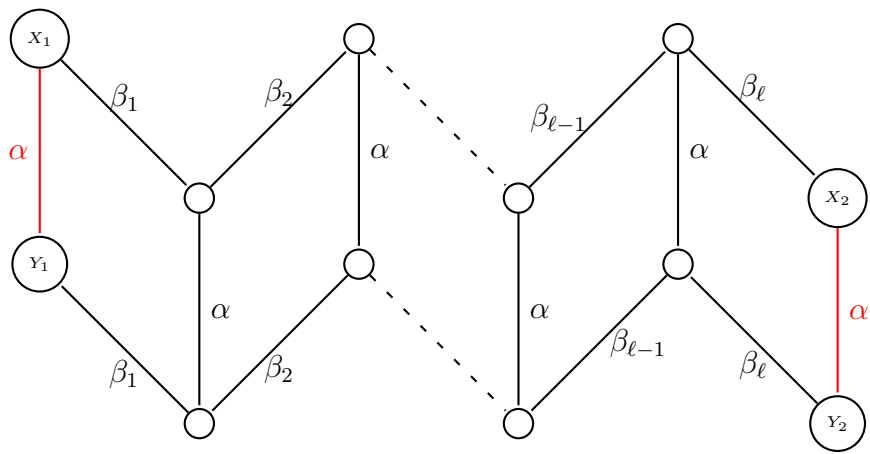

\vspace{0.2cm}
Take an element of $\mathcal{E}$ and fix a valid building
process for it. The above observations also imply, that when
observing the evolution of the inclusion graph, after the first
occurrence of an edge with some fixed label $\alpha$, new edges
with the same label can come up only when using Step $B$ always
with a different label next to $\alpha$. By easy induction on the
number of building steps, this results that between any two edges
with the same label $\alpha$ there is a ``path of $4$-cycles``.
See Figure \ref{fig path 4-cycles}. Note that as in Figure
\ref{fig path 4-cycles}, all the $\beta_i$'s must be different.
Along this path of $4$-cycles we also obtain a shortest path
between $X_1$ and $X_2$, and similarly between $Y_1$ and $Y_2$.

\section{Main results}\label{section main results}

The first of the main results of this paper is that the set
systems in $\mathcal{E}$, described in the previous section, are
actually all the extremal set systems of $VC$-dimension at most
$2$ and containing $\emptyset$.

\begin{thm}\label{main theorem 1}
$\mathcal{F}\subseteq 2^{[n]}$ is an extremal set system with
$dim_{VC}(\mathcal{F})\leq 2$ and $\emptyset\in \mathcal{F}$ iff
$\mathcal{F}\in \mathcal{E}$.
\end{thm}

Before turning to the proof of Theorem \ref{main theorem 1}, we
first prove a lemma about the building processes in $\mathcal{E}$,
that will play a key role further on.

\begin{lemma}\label{kiegeszit}
Suppose that $\mathcal{F}',\mathcal{F}$ are elements of
$\mathcal{E}$ such that $\mathcal{F}'\subseteq \mathcal{F}$. Then
$\mathcal{F}'$ can be extended with valid building process to
build up $\mathcal{F}$.
\end{lemma}
\textbf{Proof:} Suppose this is not the case, and consider a
counterexample. Without loss of generality we may suppose that the
counterexample is such that $\mathcal{F}'$ cannot be continued
with any valid step towards $\mathcal{F}$. $\mathcal{F}'$ and
$\mathcal{F}$ are both extremal and so $\mathbb{G}_{\mathcal{F}'}$
and $\mathbb{G}_{\mathcal{F}}$ are both isometrically embedded, in
particular connected, hence the neighborhood of
$\mathbb{G}_{\mathcal{F}'}$ inside $\mathbb{G}_{\mathcal{F}}$ is
nonempty. Now take a closer look at the edges on the boundary of
$\mathbb{G}_{\mathcal{F}'}$.

If there would be an edge going out from
$\mathbb{G}_{\mathcal{F}'}$ with a label $\alpha \in
supp(\mathcal{F})\backslash supp(\mathcal{F}')$, then Step $A$
would apply with this label $\alpha$. On the other hand  there
cannot be an edge going into $\mathbb{G}_{\mathcal{F}'}$ with a
label $\alpha \notin supp(\mathcal{F}')$, otherwise the endpoint
of this edge inside $\mathbb{G}_{\mathcal{F}'}$ would contain
$\alpha$, what would be a contradiction.

We can therefore assume that the label of any edge on the boundary
of $\mathbb{G}_{\mathcal{F}'}$, independently of the direction of
the edge, is an element of $supp(\mathcal{F}')$. However as
$\emptyset \in \mathcal{F}'$ and $\mathbb{G}_{\mathcal{F}'}$ is
isometrically embedded, an element belongs to $supp(\mathcal{F}')$
only if it appears as an edge label in
$\mathbb{G}_{\mathcal{F}'}$. Now take an edge $(W,V)$ on the
boundary of $\mathbb{G}_{\mathcal{F}'}$ with $W\in \mathcal{F}'$,
$V\in \mathcal{F}\backslash \mathcal{F}'$ and with some label
$\alpha$, together with an edge $(X,Y)$ with the same label inside
$\mathbb{G}_{\mathcal{F}'}$. Denote the distance of the edges
$(W,V)$ and $(X,Y)$ by $\ell$, i.e.
$d_{\mathbb{H}_n}(W,X)=d_{\mathbb{H}_n}(V,Y)=\ell$. The latter
equality means, that depending on the direction of the edges, $W$
and $X$ both do contain the element $\alpha$, or neither of them
does. Suppose that the triple $\alpha$, $(W,V)$, $(X,Y)$ is such
that the distance $\ell$ is minimal.

First suppose that $\ell>1$. Since the edges $(W,V),(X,Y)$ have
the same label and $\mathcal{F}\in \mathcal{E}$, there is a path
of $4$-cycles of length $\ell$ between them inside
$\mathbb{G}_{\mathcal{F}}$. This path of $4$-cycles also provides
shortest paths between the endpoints of the edges $(W,V),(X,Y)$.
By the minimality of our choice, in this path, except the edges at
the ends, there cannot be an edge with label $\alpha$ neither
totally inside $\mathbb{G}_{\mathcal{F}'}$, neither on the
boundary of it, meaning that this path of $4$-cycles is
essentially going outside $\mathbb{G}_{\mathcal{F}'}$. See Figure
\ref{fig l>1}.

\begin{figure}
\centering
\begin{tikzpicture}[-,>=stealth',shorten >=1pt,auto,node distance=3cm,
      thick,main node/.style={circle,fill=none,draw,font=\sffamily\tiny\bfseries}]

      \node[main node] (1) {$V$};
      \node[main node] (2) [below of=1] {$W$};
      \node[main node] (3) [above right of=1] {$P_1$};
      \node[main node] (4) [below of=3] {$Q_1$};
      \node[main node] (5) [right=3.5cm of 3] {$P_{\ell-1}$};
      \node[main node] (6) [below of=5] {$Q_{\ell-1}$};
      \node[main node] (7) [below right=4.8cm and 1.6cm of 5] {$Y$};
      \node[main node] (8) [below of=7] {$X$};
      \node[main node] (9) [left of=8] {$S$};
      \node[main node] (10) [left of=9] {$T$};

\draw[loosely dashed,-] (-2,-2)--(10,-2)
node[pos=.05,above]{$\mathbb{G}_{\mathcal{F}}\backslash
\mathbb{G}_{\mathcal{F}'}$};

\draw[loosely dashed,-] (10,-2)--(-2,-2)
node[pos=.95,below]{$\mathbb{G}_{\mathcal{F}'}$};

     \draw[] (1) -- (2) node[pos=.5,right] {$\alpha$};
     \draw[] (3) -- (4) node[pos=.5,right] {$\alpha$};
     \draw[] (5) -- (6) node[pos=.5,right] {$\alpha$};
     \draw[] (7) -- (8) node[pos=.5,right] {$\alpha$};
     \draw[loosely dashed] (3) -- (5) node[pos=.5,above] {$\beta_2,\dots,\beta_{\ell-1}$};
     \draw[loosely dashed] (4) -- (6) node[pos=.5,above] {$\beta_2,\dots,\beta_{\ell-1}$};
     \draw[] (1) -- (3) node[pos=.6,below] {$\beta_1$};
     \draw[] (2) -- (4) node[pos=.8,below] {$\beta_1$};
     \draw[] (5) -- (7) node[pos=.4,right] {$\beta_{\ell}$};
     \draw[] (6) -- (8) node[pos=.4,right] {$\beta_{\ell}$};
     \draw[] (10) -- (9) node[pos=.5,above] {$\beta_1$};
     \draw[snake] (2) -- (10) node[pos=.5,right] {};
     \draw[snake] (9) -- (8) node[pos=.5,right] {};

    \end{tikzpicture}
    \caption{Case $l>1$}
    \label{fig l>1}
\end{figure}

Since $\mathbb{G}_{\mathcal{F}'}$ is isometrically embedded and
$d_{\mathbb{H}_n}(W,X)=\ell$, there must be a path of length
$\ell$ between $W$ and $X$ inside $\mathbb{G}_{\mathcal{F}'}$. As
this path runs inside $\mathbb{G}_{\mathcal{F}'}$, it has to be
disjoint from the path of $4$-cycles. Along the path of $4$-cycles
all the $\beta_i$'s are different, so for each $i$ exactly one of
the sets $W$ and $X$ contains the element $\beta_i$. In particular
for $i=1$, the shortest path between $W$ and $X$ inside
$\mathbb{G}_{\mathcal{F}'}$ also has to contain an edge $(T,S)$
with label $\beta_1$ with direction determined by the sets $W$ and
$X$. However the distance between $W$ and $T$ is at most $\ell-1$,
and hence the triple $\beta_1,(W,Q_1),(T,S)$ contradicts with the
minimality of the initial triple $\alpha,(W,V),(X,Y)$ where the
distance was $\ell$.

By the above reasoning only $\ell=1$ is possible. In this case the
endpoints of the edges $(W,V)$, $(X,Y)$ are connected by edges
with the same label. Let this label be $\beta$. The direction of
these edges is predetermined by $\mathbb{G}_{\mathcal{F}'}$.
$\{\alpha,\beta\}\notin st(\mathcal{F}')$, otherwise there would
be already a copy of $\mathbb{G}_{2^{\{\alpha,\beta\}}}$ in
$\mathbb{G}_{\mathcal{F}'}$, which together with the vertices
$W,V,X,Y$ would give us two different copies of it inside
$\mathbb{G}_{\mathcal{F}}$, which is impossible by Observation
\ref{max el uniq st}, as $\{\alpha,\beta\}$ is a maximal set
strongly shattered by the extremal family $\mathcal{F}$. Hence
Step $B$ applies with new vertex $V$, edges $(W,V),(V,Y)$ and
labels $\alpha,\beta$ respectively, contradicting with the fact,
that we started with a counterexample. See Figure \ref{fig l=1}.
$\blacksquare$ \\

\begin{figure}
\centering
\begin{tikzpicture}[-,>=stealth',shorten >=1pt,auto,node distance=2cm,
      thick,main node/.style={circle,fill=none,draw,font=\sffamily\tiny\bfseries}]

      \node[main node] (1) {$V$};
      \node[main node] (2) [below of=1] {$W$};
      \node[main node] (3) [right of=2] {$Y$};
      \node[main node] (4) [below of=3] {$X$};

\draw[loosely dashed,-] (-1.3,-1.3)--(3.3,-1.3)
node[pos=.05,above]{$\mathbb{G}_{\mathcal{F}}\backslash
\mathbb{G}_{\mathcal{F}'}$};

\draw[loosely dashed,-] (3.3,-1.3)--(-1.3,-1.3)
node[pos=.95,below]{$\mathbb{G}_{\mathcal{F}'}$};

     \draw[] (1) -- (2) node[pos=.5,right] {$\alpha$};
     \draw[] (3) -- (4) node[pos=.5,right] {$\alpha$};
     \draw[] (1) -- (3) node[pos=.5,right] {$\beta$};
     \draw[] (2) -- (4) node[pos=.5,right] {$\beta$};

    \end{tikzpicture}
    \caption{Case $\ell=1$}
    \label{fig l=1}
\end{figure}

Now we are ready to prove Theorem \ref{main theorem 1}.\\

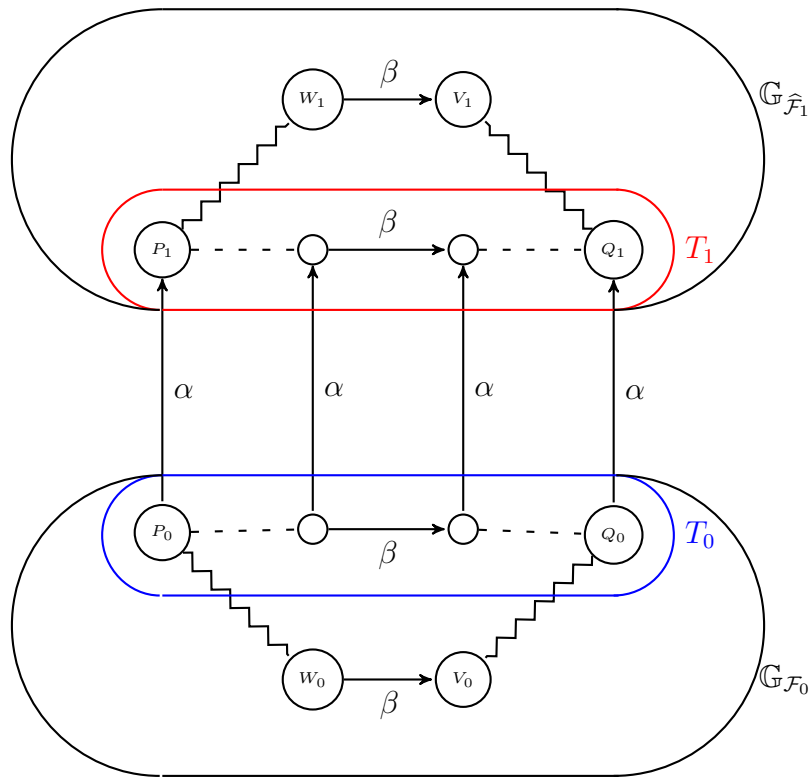
\begin{figure}
\centering
\begin{tikzpicture}[-,>=stealth',shorten >=1pt,auto,node distance=2cm,
      thick,main node/.style={circle,fill=none,draw,font=\sffamily\tiny\bfseries}]

      \node[main node] (1) {$P_1$};
      \node[main node] (2) [right of=1] {};
      \node[main node] (3) [right of=2] {};
      \node[main node] (4) [right of=3] {$Q_1$};
      \node[main node] (5) [below=3cm of 1] {$P_0$};
      \node[main node] (6) [below=3.3cm of 2] {};
      \node[main node] (7) [below=3.3cm of 3] {};
      \node[main node] (8) [below=3cm of 4] {$Q_0$};
      \node[main node] (9) [above of=2] {$W_1$};
      \node[main node] (10) [above of=3] {$V_1$};
      \node[main node] (11) [below of=6] {$W_0$};
      \node[main node] (12) [below of=7] {$V_0$};

     \draw[->] (9) -- (10) node[pos=.5,above] {$\beta$};
     \draw[->] (11) -- (12) node[pos=.5,below] {$\beta$};
     \draw[<-] (1) -- (5) node[pos=.5,right] {$\alpha$};
     \draw[<-] (2) -- (6) node[pos=.5,right] {$\alpha$};
     \draw[<-] (3) -- (7) node[pos=.5,right] {$\alpha$};
     \draw[<-] (4) -- (8) node[pos=.5,right] {$\alpha$};
     \draw[->] (2) -- (3) node[pos=.5,above] {$\beta$};
     \draw[->] (6) -- (7) node[pos=.5,below] {$\beta$};
     \draw[snake] (1) -- (9) node[pos=.6,right] {};
     \draw[snake] (4) -- (10) node[pos=.8,right] {};
     \draw[snake] (5) -- (11) node[pos=.4,right] {};
     \draw[snake] (8) -- (12) node[pos=.4,right] {};
     \draw[loosely dashed] (1) -- (2) node[pos=.6,right] {};
     \draw[loosely dashed] (3) -- (4) node[pos=.8,right] {};
     \draw[loosely dashed] (5) -- (6) node[pos=.4,right] {};
     \draw[loosely dashed] (7) -- (8) node[pos=.4,right] {};

    \draw[red] (0,0.8) arc (90:270:0.8);
    \draw[red] (6,-0.8) arc (-90:90:0.8);
    \draw[red] (0,0.8) -- (6.1,0.8);
    \draw[red] (0,-0.8) -- (6.1,-0.8);

    \draw[blue] (0,-3.0) arc (90:270:0.8);
    \draw[blue] (6,-4.6) arc (-90:90:0.8);
    \draw[blue] (0,-3.0) -- (6.1,-3.0);
    \draw[blue] (0,-4.6) -- (6.1,-4.6);

    \draw[] (0,3.2) arc (90:270:2);
    \draw[] (6,-0.8) arc (-90:90:2);
    \draw[] (0,3.2) -- (6.1,3.2);

    \draw[] (0,-3.0) arc (90:270:2);
    \draw[] (6,-7.0) arc (-90:90:2);
    \draw[] (0,-7.0) -- (6.1,-7.0);

    \node[red] [right=0.4 of 4]{$T_1$};
    \node[blue] [right=0.4 of 8]{$T_0$};
    \node[] [right=3.4 of 10]{$\mathbb{G}_{\widehat{\mathcal{F}}_1}$};
    \node[] [right=3.4 of 12]{$\mathbb{G}_{\mathcal{F}_0}$};

    \end{tikzpicture}
    \caption{Building up extremal set systems}
    \label{build extr syst}
\end{figure}

\textbf{Proof:} One direction of the theorem is just Lemma
\ref{steps keep extremality}. For the other direction we use
induction on the number of sets in $\mathcal{F}$. If
$|\mathcal{F}|=1$, then $\mathcal{F}$ is necessarily
$\{\emptyset\}$, and so belongs trivially to $\mathcal{E}$. Now
suppose we proved the statement for all set systems with at most
$m-1$ members, and let $\mathcal{F}$ be an extremal family of size
$m$, of $VC$-dimension at most $2$ and containing $\emptyset$.
Take an arbitrary element $\alpha$ appearing as a label of an edge
going out from $\emptyset$ in $\mathbb{G}_{\mathcal{F}}$, i.e. an
element $\alpha$ such that $\{\alpha\}\in \mathcal{F}$. Consider
the standard subdivision of $\mathcal{F}$ with respect to the
element $\alpha$ with parts $\mathcal{F}_0$ and $\mathcal{F}_1$
(see Definition \ref{stand subdiv def}), and let
\[\widehat{\mathcal{F}}_1=\{F\cup\{\alpha\}\ :\ F\in
\mathcal{F}_1\}.\] Note that with respect to shattering and strong
shattering $\mathcal{F}_1$ and $\widehat{\mathcal{F}}_1$ behave in
the same way. Since $\mathcal{F}$ is extremal, so are
$\mathcal{F}_0$, $\mathcal{F}_1$ and hence
$\widehat{\mathcal{F}}_1$ as well, and clearly their
$VC$-dimension is at most $2$. The collection of all edges with
label $\alpha$ in the inclusion graph $\mathbb{G}_{\mathcal{F}}$
forms a cut. This cut divides $\mathbb{G}_{\mathcal{F}}$ into two
parts, that are actually the inclusion graphs
$\mathbb{G}_{\mathcal{F}_0}$ and
$\mathbb{G}_{\widehat{\mathcal{F}}_1}$. Note that
$\mathbb{G}_{\mathcal{F}_1}$ and
$\mathbb{G}_{\widehat{\mathcal{F}}_1}$ are isomorphic as directed
edge labelled graphs. Let $T_0$ and $T_1$ be the induced subgraphs
on the endpoints of the cut edges in $\mathbb{G}_{\mathcal{F}_0}$
and $\mathbb{G}_{\widehat{\mathcal{F}}_1}$, respectively. See
Figure \ref{build extr syst}. $T_0$ and $T_1$ are isomorphic, and
they are actually the inclusion graphs of the set systems
$\mathcal{T}_0=\mathcal{F}_0\cap \mathcal{F}_1$ and
$\mathcal{T}_1=\{F\cup\{\alpha\}, F\in \mathcal{T}_0\}$. Similarly
to the pair $\mathcal{F}_1$, $\widehat{\mathcal{F}}_1$, the set
systems $\mathcal{T}_0$ and $\mathcal{T}_1$ also behave in the
same way with respect to shattering and strong shattering. By
assumption $\mathcal{F}$ is extremal, and so according to
Proposition 5.1 from \cite{VC1} so is $\mathcal{T}_0$ and hence
$\mathcal{T}_1$. For every set $S$ in
$Sh(\mathcal{T}_0)=Sh(\mathcal{F}_0\cap \mathcal{F}_1)\subseteq
2^{[n]\backslash \{\alpha\}}$ the set $S\cup\{\alpha\}$ is
shattered by $\mathcal{F}$, implying that
$dim_{VC}(\mathcal{T}_0)\leq dim_{VC}(\mathcal{F})-1\leq 1$.
Therefore $\mathcal{T}_0$ is an extremal family of $VC$-dimension
at most $1$, and so by Proposition \ref{extr char VC-dim=1} we get
that $T_0$ (and hence $T_1$) is a directed edge labelled tree
having all edge labels different. Note that for any edge label
$\beta$ appearing in $T_0$ (and hence in $T_1$), there is a copy
of $\mathbb{G}_{2^{\{\alpha,\beta\}}}$ along the cut, implying
that $\{\alpha,\beta\}\in st(\mathcal{F})=Sh(\mathcal{F})$. By the
$VC$-dimension constraint on $\mathcal{F}$ the set
$\{\alpha,\beta\}$ is a maximal element of
$st(\mathcal{F})=Sh(\mathcal{F})$, and so by Observation \ref{max
el uniq st} there cannot be another copy of
$\mathbb{G}_{2^{\{\alpha,\beta\}}}$ in $\mathbb{G}_{\mathcal{F}}$,
neither in $\mathbb{G}_{\mathcal{F}_0}$ nor in
$\mathbb{G}_{\widehat{\mathcal{F}}_1}$, in particular
$\{\alpha,\beta\}\notin st(\mathcal{F}_0)$.

Let's now turn to the building process of $\mathcal{F}$. Our
choice of $\alpha$ guarantees that $\emptyset \in \mathcal{F}_0,
\mathcal{F}_1$ and so by the induction hypothesis both of them
belong to $\mathcal{E}$. In particular we can build up
$\mathcal{F}_0$, and in the meantime $\mathbb{G}_{\mathcal{F}_0}$,
according to the building rules in $\mathcal{E}$. $\alpha\notin
supp(\mathcal{F}_0)$ and so we can apply Step $A$ with $\alpha$ to
add one fixed cut edge to $\mathbb{G}_{\mathcal{F}_0}$. Then we
apply Step $B$ several times to add the whole of $T_1$ to
$\mathbb{G}_{\mathcal{F}_0}$ and simultaneously $\mathcal{T}_1$ to
$\mathcal{F}_0$. By earlier observations all edge labels of $T_1$
are different, and if $\beta$ is such a label, then
$\{\alpha,\beta\}\notin st(\mathcal{F}_0)$, and hence all these
applications of Step $B$ will be valid ones. The building process
so far shows that $\mathcal{F}_0\cup\mathcal{T}_1$ is also a
member of $\mathcal{E}$.
$\mathbb{G}_{\mathcal{F}_0\cup\mathcal{T}_1}$ is just
$\mathbb{G}_{\mathcal{F}_0}$ and $T_1$ glued together along the
cut in the way described above.

$T_0$ shows that $\mathcal{T}_0$ can be built up using only Step
$A$, and hence it belongs to $\mathcal{E}$. The inclusion
$\mathcal{T}_1\subseteq \widehat{\mathcal{F}}_1$ shows that
$\mathcal{T}_0\subseteq \mathcal{F}_1$, therefore by Lemma
\ref{kiegeszit} $\mathcal{T}_0$ can be extended with a valid
building process to build up $\mathcal{F}_1$. This extension can
also be considered as building up $\widehat{\mathcal{F}}_1$ from
$\mathcal{T}_1$. $\emptyset\notin \mathcal{T}_1,
\widehat{\mathcal{F}}_1$ and so neither of the two systems is a
member of $\mathcal{E}$, however this causes no problems, as the
pairs $\mathcal{T}_0$, $\mathcal{T}_1$ and $\mathcal{F}_1$,
$\widehat{\mathcal{F}}_1$ behave in the same way with respect to
shattering and strong shattering, and so all building steps remain
valid.

We claim, that this last building procedure remains valid, and so
completes a desired building process for $\mathcal{F}$, if we
start from $\mathcal{F}_0\cup\mathcal{T}_1$ instead of
$\mathcal{T}_1$. First note that if there is a label appearing
both in $\mathbb{G}_{\mathcal{F}_0}$ and
$\mathbb{G}_{\widehat{\mathcal{F}}_1}$, then it appears also in
$T_0$, and hence in $T_1$. Indeed let $\beta$ be such a label, and
consider $2$ edges with this label, one going from $W_0$ to $V_0$
in $\mathbb{G}_{\mathcal{F}_0}$ and the other going from $W_1$ to
$V_1$ in $\mathbb{G}_{\widehat{\mathcal{F}}_1}$. See Figure
\ref{build extr syst}. $\mathbb{G}_{\mathcal{F}}$ is isometrically
embedded, therefor there is a shortest path both between $W_0$ and
$W_1$ and between $V_0$ and $V_1$ in $\mathbb{G}_{\mathcal{F}}$.
Thanks to $\beta$ these two paths have to be disjoint. Both of
these paths must have a common edge with the cut, say $(P_0,P_1)$
and $(Q_0,Q_1)$, with $P_0$ and $Q_0$ in
$\mathbb{G}_{\mathcal{F}_0}$. Since $\beta\in P_0\bigtriangleup
Q_0$, along the shortest path between $P_0$ and $Q_0$ in the
isometrically embedded inclusion graph $T_0$ of the extremal
family $\mathcal{T}_0$ there must be an edge with label $\beta$.
According to this, when applying Step $A$ in the extension
process, then the used element will be new not just when we start
from $\mathcal{T}_1$, but also when starting from
$\mathcal{F}_0\cup\mathcal{T}_1$.

Finally suppose that an application of Step $B$ with some labels
$\beta,\gamma$ in the extension process turns invalid when we
start from $\mathcal{F}_0\cup\mathcal{T}_1$ instead of
$\mathcal{T}_1$. This is possible only if $\{\beta,\gamma\}\in
st(\mathcal{F}_0\cup \mathcal{T}_1)\backslash st(\mathcal{T}_0)$,
i.e. there is a copy of $\mathbb{G}_{2^{\{\beta,\gamma\}}}$
already in $\mathbb{G}_{\mathcal{F}_0\cup\mathcal{T}_1}$. However
this copy together with the copy, that the invalid use of Step $B$
results, gives two different occurrences of
$\mathbb{G}_{2^{\{\beta,\gamma\}}}$ inside
$\mathbb{G}_{\mathcal{F}}$, which is impossible by Observation
\ref{max el uniq st}, as $\{\beta,\gamma\}$ is a maximal set
strongly
shattered by the extremal family $\mathcal{F}$. $\blacksquare$\\

As a corollary of Theorem \ref{main theorem 1} one can solve an
open problem, posed in \cite{VC1}, in the special case when the
$VC$-dimension of the systems investigated is bounded by $2$.

\begin{conj}\label{bovit sejt}
(See \cite{VC1}) For a nonempty s-extremal family
$\mathcal{F}\subseteq 2^{[n]}$ does there always exist a set $F\in
\mathcal{F}$ such that $\mathcal{F}\backslash \{F\}$ is still
s-extremal?
\end{conj}

The case when the $VC$-dimension of the systems investigated is
bounded by $1$ was solved in \cite{VC1}. Here we propose a
solution for set systems of $VC$-dimension at most $2$.

\begin{thm}\label{main theorem 2}
Let $\mathcal{F}\subseteq 2^{[n]}$ be a nonempty extremal family
of $VC$ dimension at most $2$. Then there exists an element $F\in
\mathcal{F}$ such that $\mathcal{F}\backslash\{F\}$ is still
extremal.
\end{thm}
\textbf{Proof:} Let $F\in \mathcal{F}$ be an arbitrary set from
the set system. Recall that $\varphi_i$ is the $i$th bit flip
operation, and let $\varphi=\prod_{i\in F} \varphi_i$. Since bit
flips preserve extremality, $\varphi(\mathcal{F})$ is extremal as
well. Moreover $\varphi(F)=\emptyset\in \varphi(\mathcal{F})$, and
so by Theorem \ref{main theorem 1} we have
$\varphi(\mathcal{F})\in \mathcal{E}$, hence we can consider a
building process for it. Let $V\in \varphi(\mathcal{F})$ be the
set added in the last step of this building process. The same
building process shows that
$\varphi(\mathcal{F})\backslash\{V\}\in\mathcal{E}$, and hence by
Theorem \ref{main theorem 1} we have that
$\varphi(\mathcal{F})\backslash\{V\}$ is an extremal family of
$VC$ dimension at most $2$ and containing $\emptyset$. However
$\varphi(\mathcal{F})\backslash\{V\}=\varphi(\mathcal{F}\backslash\{\varphi(V)\})$,
and since bit flips preserve extremality, we get that
$\varphi(\varphi(\mathcal{F}\backslash\{\varphi(V)\}))=\mathcal{F}\backslash\{\varphi(V)\}$
is also extremal, meaning that the set $\varphi(V)\in \mathcal{F}$
can be removed from the extremal system $\mathcal{F}$ so that the
result is still extremal. $\blacksquare$\\

\section{Concluding remarks and future work}\label{section remarks}

The building process from Section \ref{preliminaries} can be
generalized to the case when the $VC$-dimension bound is some
fixed natural number $t$. We can define a building step for every
set $S\subseteq [n]$ with $|S|\leq t$. Let Step$(\emptyset)$ be
the initialization, after which we are given the set system
$\{\emptyset\}$. For some set $S\subseteq [n]$ with $|S|\leq t$,
Step($S$) can be applied to a set system $\mathcal{F}$, if there
exists some set $F\subseteq [n]$, $F\notin \mathcal{F}$, such that
$S\in st(\mathcal{F}\cup\{F\})\backslash st(\mathcal{F})$. If such
set $F$ exists, choose one, and let the resulting system be
$\mathcal{F}\cup\{F\}$. In terms of the inclusion graph $S\in
st(\mathcal{F}\cup\{F\})\backslash st(\mathcal{F})$ means, that by
adding the set $F$ there arises a copy of $\mathbb{G}_{2^S}$
inside $\mathbb{G}_{\mathcal{\mathcal{F}\cup\{F\}}}$ containing
the vertex $F$. Similarly as previously, one can prove that $F$'s
only neighbors are the ones contained in this copy of
$\mathbb{G}_{2^S}$. Using this observation Step($S$) could have
been defined in terms of the inclusion graph as well (as it was
done in the case $t=2$).

Restrict our attention to those set systems, that can be built up
starting with Step($\emptyset$), and then using always new
building steps, i.e. not using a building step with the same set
$S$ twice. Along the same lines of thinking as in the case $t=2$,
one can prove that every such set system is extremal. We think,
that these set systems are actually all the extremal families of
$VC$-dimension at most $t$. Unfortunately, for the time being we
were unable to prove a suitable generalization of Lemma
\ref{kiegeszit}. Once it is done, the generalization of Theorem
\ref{main theorem 1}, and as a corollary a generalization of
Theorem \ref{main theorem 2} would follow easily. Although the
general version Theorem \ref{main theorem 1} would not give such a
transparent structural description of extremal systems as in the
case $t=1$, but still, its corollary, the generalization of
Theorem \ref{main theorem 2} would solve the open problem proposed
in \cite{VC1} in its entire generality.

\subsection*{Acknowledgements}
We thank L\'aszl\'o Kozma and Shay Moran for pointing out an error
in an earlier version a of the manuscript.

\end{document}